\documentclass[reqno, 11pt]{article}
\usepackage[centertags,tbtags]{amsmath}
\usepackage{amssymb} 
\usepackage{amsthm}
\usepackage{txfonts}
\usepackage{fullpage}
\usepackage[usenames,dvipsnames]{color}

\usepackage[pdftex]{thumbpdf}       
\usepackage{ifpdf}
\ifpdf
\usepackage[pdftex]{hyperref}
\hypersetup{%
pdftitle={Renewal Processes with Costs and Rewards},
pdfauthor={M. Vlasiou},
pdfkeywords={EORMS},
bookmarks=true,%
bookmarksnumbered=true,
pdfstartview={FitH},
colorlinks=true,
citecolor=Plum,
linkcolor=RedOrange,
urlcolor=RawSienna,
bookmarksopen=true,%
bookmarksopenlevel=1,     
unicode=true,%
breaklinks=true,%
}%
\else
\newcommand\phantomsection\relax
\newcommand{\url}[1]{#1}
\newcommand{\href}[2]{#2}
\usepackage{nohyperref}
\fi

\theoremstyle{plain}              
\newtheorem{theorem}{Theorem}

\theoremstyle{definition}

\newtheorem{ex}{Example}

\newcommand{\e}{\mathbb{E}}
\newcommand{\p}{\mathbb{P}}
\renewcommand{\d}{\,\mathrm{d}}

\begin{document}
\title{Renewal Processes with Costs and Rewards}
\author{Maria Vlasiou\thanks{Dept.\ of Mathematics \& Computer Science, Eindhoven University of Technology, P.O.\ Box 513, 5600 MB Eindhoven, The Netherlands, \href{mailto:m.vlasiou@tue.nl}{m.vlasiou@tue.nl}}}
\date{\today}
\maketitle

\begin{abstract}
We review the theory of renewal reward processes, which describes renewal processes that have some cost or reward associated with each cycle. We present a new simplified proof of the renewal reward theorem that mimics the proof of the elementary renewal theorem and avoids the technicalities in the proof that is presented in most textbooks. Moreover, we mention briefly the extension of the theory to partial rewards, where it is assumed that rewards are not accrued only at renewal epochs but also during the renewal cycle. For this case, we present a counterexample which indicates that the standard conditions for the renewal reward theorem are not sufficient; additional regularity assumptions are necessary. We present a few examples to indicate the usefulness of this theory, where we prove the inspection paradox and Little's law through the renewal reward theorem.
\end{abstract}

\phantomsection
\addcontentsline{toc}{section}{Basic notions and results}
\section*{Basic notions and results}

Many applications of the renewal theory involve rewards or costs (which can be simply seen as a negative reward). For example, consider the classical example of a renewal process: a machine component gets replaced upon failure or upon having operated for $T$ time units. Then the time the $n$-th component is in service is given by $Y_n=\min\{X_n,T\}$, where $X_n$ is the life of the component. In this example, one might be interested in the rate of the number of replacements in the long run. An extension of this basic setup is as follows. A component that has failed will be replaced at a cost $c_f$, while a component that is replaced while still being operational (and thus at time $T$) costs only $c<c_f$. In this case, one might be interested in choosing the optimal time $T$ that minimises the long-run operational costs. The solution to this problem involves the analysis of \emph{renewal processes with costs and rewards}.

Motivated by the above, let $\{N(t), t\geqslant 0\}$ be a renewal process with interarrival times $X_n$, $n\geqslant1$, and denote the time of the $n$-th renewal by $S_n=X_1+\cdots + X_n$. Now suppose that at the time of each renewal a reward is received; we denote by $R_n$ the reward received at the end of the $n$-th cycle. We further assume that $(R_n, X_n)$ is a sequence  of i.i.d.\ random variables, which allows for $R_n$ to depend on $X_n$. For example, $N(t)$ might count the number of rides a taxi gets up to time $t$. In this case, $X_n$ is the length of each trip and one reasonably expects the fare $R_n$ to depend on $X_n$. As usual, we denote by $(R,X)$ the generic bivariate random variable that are distributed identically to the sequence of rewards and interarrival times $(R_n,X_n)$. In the analysis, together with the standard assumption in renewal processes that the interarrival times have a finite expectation $\e[X]=\tau$, we will further assume that $\e[\,|R|\,]<\infty$. The cumulative reward up to time $t$ is given by $R(t)=\sum_{n=1}^{N(t)} R_n$, where the sum is taken to be equal to zero in the event that $N(t)=0$. Depending on whether the reward is collected at the beginning or at the end of the renewal period, one might wish to adapt this definition by taking $R(t)=\sum_{n=1}^{N(t)+1} R_n$.

The importance of renewal processes with costs and rewards is evidenced by their application to Markovian processes and models. Most queuing processes in which customers arrive according to a renewal process are regenerative processes (see the article on \textit{Regenerative Processes}) with cycles beginning each time an arrival finds the system empty. Moreover, for every regenerative process $Y(t)$ we may define a reward structure as follows: $R_n=\int_{S_{n-1}}^{S_n} Y(t) \d t$, where  all $R_n$, $n \geqslant1$, are i.i.d., except possibly for $R_1$ that might follow a different distribution. The basic tools that are used are the computation of the reward per unit of time and the rate of the expected value of the reward. These two results are summarised in the following theorem.

\begin{theorem}\label{th:1}
Let $N(t)$ be a renewal reward process generated by $(R, X)$, $n\geqslant1$. Assume that $\e[\,|R|\,]<\infty$ and let $r= \e[R]$, $\tau=\e[X]<\infty$. Then
\begin{align}
\label{eq:rv}
\lim_{t\to\infty} \frac{R(t)}{t}&= \frac{r}{\tau}, \quad \mbox{with probability 1, and}\\
\label{eq:exp}
\lim_{t\to\infty} \frac{\e[R(t)]}{t}&= \frac{r}{\tau}.
\end{align}
\end{theorem}

These results do not depend on whether rewards are collected at the beginning of the renewal cycle or at its end, or if they are collected in a more complicated fashion (e.g.\ continuously or in a non-monotone way during the cycle), as long as $(R_n, X_n)$,  $n\geqslant1$, is a sequence  of i.i.d.\ bivariate random variables and minor regularity conditions hold; see the following section. Moreover, the results remain valid also for delayed renewal reward processes, i.e.\ when the first cycle might follow a different distribution.

Note that \eqref{eq:exp} does not follow from \eqref{eq:rv} since almost sure convergence does not imply the convergence of the expected values. Think of the distribution
$$
\p[Y_n=x]=\begin{cases}
  n/(n+1), &x=0,\\
  1/(n+1), &x=n+1.
\end{cases}
$$
Since for all $n$, $\e[Y_n]$=1, we have $\lim_{n\to\infty} \e[Y_n]=1$. Moreover, $Y_n$ goes almost surely to zero, i.e.\ $\p[\lim_{n\to\infty} Y_n = 0]=1$, which means that the value of the random variable and its expectation differ in the limit.
\begin{proof}[Proof of Theorem \ref{th:1}]
To prove \eqref{eq:rv}, write
$$
\frac{R(t)}{t}=\frac{\sum_{n=1}^{N(t)} R_n}{N(t)}\frac{N(t)}{t}.
$$
By the strong law of large numbers we have that
$$
\lim_{t\to\infty} \frac{\sum_{n=1}^{N(t)} R_n}{N(t)}= \e[R]=r,
$$
since $N(t)$ goes to infinity almost surely, and by the strong law for renewal processes we know that
$$
\lim_{t\to\infty} \frac{N(t)}{t}=\frac{1}{\e[X]}=\frac{1}{\tau}.
$$

The proof of \eqref{eq:exp} is a bit more involved. The classical proof applies Wald's equation \cite{wald44} and requires the proof that $\lim_{t\to\infty}\frac{\e[R_{N(t)+1}]}{t}=0$ by constructing a renewal-type equation (see the article on \textit{Renewal Function and Renewal-Type Equations}). Here we present an alternative proof that follows the same steps of the classical proof for the Elementary Renewal Theorem (see e.g.\ \cite{asmussen-APQ}) but seems to be new in the setting of renewal reward processes.

We first construct a lower bound. From Fatou's lemma we have that
\begin{equation}
\liminf_{t\to\infty} \frac{\e[R(t)]}{t} \geqslant \e [\liminf_{t\to\infty} \frac{R(t)}{t}]=\e[\lim_{t\to\infty} \frac{R(t)}{t}]=\frac{r}{\tau}
\end{equation}
from \eqref{eq:rv}.

For the upper bound we use truncation. Take $M<\infty$ and set $R_n^M=\max\{R_n, -M\}$. Note that
$$
R(t)=\sum_{n=1}^{N(t)} R_n\leqslant\sum_{n=1}^{N(t)} R_n^M \stackrel{def}{=} R^M(t).
$$
We then have that
$$
\frac{\e[R(t)]}{t}\leqslant \frac{\e[R^M(t)]}{t}=\frac{\e[\sum_{n=1}^{N(t)+1} R_n^M - R_{N(t)+1}^M]}{t}=\frac{\e[N(t)+1]\e[R_1^M]}{t}-\frac{\e[R_{N(t)+1}^M]}{t}\leqslant \frac{\e[N(t)]+1}{t}\e[R_1^M] +\frac{M}{t}.
$$
The last equality follows from Wald's equation since $N(t)+1$ is a stopping time and for the last inequality we simply observe that we can bound $-\e[R_n^M]$ from above with $M$ (by construction). Thus,
$$
\limsup_{t\to\infty} \frac{\e[R(t)]}{t}\leqslant\limsup_{t\to\infty} \left( \frac{\e[N(t)]+1}{t}\e[R_1^M] +\frac{M}{t}\right)=\frac{\e[R_1^M]}{\tau}
$$
by the Elementary Renewal Theorem. Since we can choose any $M<\infty$ in our construction, by monotone convergence we have that $\e[R_1^M]\to\e[R_1]=r$ as $M\to\infty$.

The advantage of this proof is that it avoids the computation of $\lim_{t\to\infty}\frac{\e[R_{N(t)+1}]}{t}$ and that it is in essentials identical to the classical proof for the corresponding theorem in renewal processes without rewards.
\end{proof}

\section*{Partial Rewards}
\phantomsection
\addcontentsline{toc}{section}{Partial Rewards}

In some applications rewards might not be given at the beginning or an end of a cycle, but might be earned in a continuous (maybe non-monotone) fashion during a cycle. The simplest case is when the reward during a cycle is given at a constant rate during the cycle, or only a portion thereof. The portion of the reward up to time $t$ associated with a cycle starting at $X_{N(t)}$ is called a \emph{partial reward}. Recall that we have assumed that $E[X]<\infty$ and $\e[\,|R|\,]<\infty$. Under further assumptions, Theorem~\ref{th:1} still holds. In the case of partial rewards, we need to assume some form of good behaviour as nothing so far prevents the rewards from escaping to plus or minus infinity while still having the average reward per cycle being finite. As an example, to avoid such erratic behaviour consider the restrictive assumption that the absolute value of rewards is uniformly bounded by some number; then Theorem~\ref{th:1} holds. This is also the case when rewards are non-negative or if they accumulate in a monotone fashion over any interval. Wolff \cite{wolff-SMTQ} assumes that the partial reward during a cycle is the difference of two different monotone non-negative processes, which could be interpreted as income and costs over the renewal cycle.

It seems to be a standard misconception that only having $\e[\,|R|\,]<\infty$ is sufficient for Theorem~\ref{th:1} to hold also for partial rewards. In lack of necessary and sufficient conditions in the literature for for this result, or at least of a counterexample where $\e[\,|R|\,]<\infty$ but where Theorem~\ref{th:1} does not hold for partial rewards, we present here a counterexample that illustrates the needs for additional assumptions. In this case, the average reward over a cycle is finite (and we will construct it to be equal to zero), but the average reward for every finite $t$ is not defined, and thus also not its limit over $t$; see \eqref{eq:exp}.

Take a renewal reward process that is defined as follows: $X_i=1$ and $R_i=0$ for all $i$. During a cycle, rewards build up and are depleted in the following fashion. For cycle $j$ choose a Cauchy random variable $C_j$ and take the auxiliary function $f(x)=x$ for $x\in[0,1/2]$ and $f(x)=1-x$ for $x\in[1/2,1]$. Define the reward for any time $t$ during cycle $j$ as $C(t)=C_j \cdot f(t-[t])$. Then we see that the cumulative reward up to time $t$, for all $t>0$ is given by $R(t)=C(t)$. Thus, since a Cauchy random variable has not a well-defined first moment, we see that the average cumulative reward in this example cannot be defined, and thus Theorem~\ref{th:1} cannot be extended for this case, despite the fact that at the end of a cycle (and thus at all integer times) the cumulative reward is equal to zero.

In the following, we will see two examples where partial rewards are assumed.

\section*{Examples}
\phantomsection
\addcontentsline{toc}{section}{Examples}
We proceed with a few applications of renewal processes with costs and rewards that exhibit the usefulness of these processes. The first two show how (cyclic) renewal processes can be seen as a special case of renewal processes with costs and rewards, while the last two make use of rewards earned continuously over time. Naturally, an abundance of examples of renewal (reward) processes might be found in Markov processes and chains.

\begin{ex}[Standard renewal processes]
Observe that the standard renewal process is simply a special case of a renewal reward process, where we assume that the reward at each cycle is equal to 1. In this case we see for example that \eqref{eq:exp} is simply the Elementary Renewal Theorem.
\end{ex}

\begin{ex}[Alternating renewal processes]
For an alternating renewal process (see the article on \textit{Alternating Renewal Processes}), i.e.\ a renewal process that can be in one of two phases (say ON and OFF) suppose that we earn at the end of a cycle an amount equal to the time the system was ON during that cycle. Alternatively, one may assume that we earn at a rate of one per unit of time when the system is ON, but one need not consider continuous rewards at this point. Then the total reward earned in the interval $[0,t]$ is equal to the total ON time in that interval, and thus by \eqref{eq:rv} we have that as $t\to\infty$
$$
\frac{\mbox{ON time in} [0,t]}{t}\to\frac{\e[Y_{\rm ON}]}{\e[Y_{\rm ON}]+\e[Y_{\rm OFF}]},
$$
where $Y_{\rm ON}$ is a generic ON time in a cycle and equivalently for the OFF times. Thus, for non-lattice distributions, the limiting probability of the system being ON is equal to the long-run proportion of time it is ON. Naturally, these results extend to general cyclic renewal processes, i.e.\ a process with more than two phases \cite{serfozo-BASP}.
\end{ex}

\begin{ex}[Average time of age and excess]
Let $A(t)$ denote the age at time $t$ of a renewal process generated by $X_n \sim X$ and suppose that we are interested in computing the average value of the age. With a slight abuse of terminology, define this to be equal to
$$
\lim_{t\to\infty} \frac{1}{t} \int_0^t A(s) \d s.
$$
Suppose now that we are obtaining a reward continuously at a rate equal to the age of the renewal process at that time. Thus, $\int_0^t A(s) \d s$ represents the total earnings by time $t$. Moreover, since the age of a renewal process at time $s$ since the last renewal is simply equal to $s$ we have that the reward during a renewal cycle of length $X$ is equal to $\int_0^X s \d s= X^2/2$.

Then by \eqref{eq:rv} we have that with probability 1,
$$
\lim_{t\to\infty} \frac{\int_0^t A(s) \d s}{t} = \frac{\e[X^2/2]}{\e[X]}.
$$
We can apply the same logic to the residual life at time $t$, $B(t)$, and we will end up at the same result. Since the renewal epoch covering time $t$ is equal to $X_{N(t)+1} = A(t)+B(t)$ we have that
$$
\lim_{t\to\infty} \frac{1}{t} \int_0^t X_{N(t)+1} \d s= \frac{\e[X^2]}{\e[X]}\geqslant \e[X]
$$
where we have equality only when Var$(X)=0$. Thus, we retrieve the inspection paradox.
\end{ex}

\begin{ex}[Little's law]
Consider a GI/GI/1 stable queue with the interarrival times $X_i$ generating a nonlattice renewal process. Suppose that we start observing the system upon the arrival of a customer. Denote a generic interarrival time by $X$ and let $\e[X]=1/\lambda $. Also denote by $n(t)$ the number of customers in the system at time $t$. Consider the average number of customers in the system in the long run, which we denote by $L$. Then, with a slight abuse of notation we have that
$$
L=\lim_{t\to\infty}\frac1t\int_0^t n(y) \d y.
$$
Define a renewal reward process with cycles $C$ starting each time an arrival finds the system empty and suppose that we earn a reward at time $y$ with a rate $n(y)$. Then from \eqref{eq:exp} we have that
\begin{equation}\label{eq:L}
L=\frac{\e[\mbox{reward during a cycle}]}{\e[\mbox{cycle time}]}=\frac{\e[\int_0^C n(y) \d y]}{\e[C]}.
\end{equation}

Now let $N$ be the number of customers served during a cycle and define the long-run average sojourn time of customers as $T=\lim_{n\to\infty} (T_1+\cdots+T_n)/n$, where $T_i$ is the time customer $i$ spent in the system. Suppose now that in each cycle of length $C$ where $N$ customers were served. We can see then that the reward we have received in that cycle, which was defined to be equal to a rate of $n(y)$ at each time $y$, can also be seen as having each customer pay at a rate 1 for each time unit he is in the system. In other words, in a cycle with $N$ customers, the reward is equal to $T_1+\cdots+T_N$, i.e.\ the total time all customers in that cycle spent in the system. Then $T$ is the average reward per unit time (observe that since the second definition of the reward depends on the number of customers in a cycle, then the duration of the cycle is also defined in terms of customers) and again by \eqref{eq:exp} we have that
\begin{equation}\label{eq:S}
S=\frac{\e[\mbox{reward during a cycle}]}{\e[\mbox{cycle time}]}=\frac{\e[\sum_1^N T_i]}{\e[N]}.
\end{equation}

We can now prove Little's law, one of the fundamental relations in queuing theory, which states that $L=\lambda T$. To see this, observe that
$$
C=\sum_1^N X_i,
$$
and since $N$ is a stopping time for the sequence $\{X_i\}$ we have from Wald's equation that
$$
\e[C]=\e[N]\e[X]=\e[N]/\lambda .
$$
Thus, by \eqref{eq:L} and \eqref{eq:S} we see that
$$
L=\lambda T\frac{\e[\int_0^C n(y) \d y]}{\e[\sum_1^N T_i]}.
$$
However, as we have seen, the fraction is equal to 1, since both the numerator and the denominator describe the reward earned during a cycle. Thus we retrieve Little's law.
\end{ex}

\section*{Further reading}
\phantomsection
\addcontentsline{toc}{section}{Further reading}
Almost all books on stochastic processes and introductory probability have a section on renewal theory. Here we only refer to the ones that have a separate mention of renewal processes with costs and rewards. Cox \cite{cox-RT} is one of the few manuscripts exclusively devoted to renewal processes. There, renewal processes with costs and rewards are found under the term \textit{cumulative processes}. Another term that has been used to describe these processes is \textit{compound renewal processes}. The general concept of a cumulative process and the asymptotic results related to them
are due to Smith \cite{smith55,smith58}. Mercer and Smith \cite{mercer59} investigate cumulative processes acclimated with a Poisson process, in connection with a study of the wear of conveyor
belting.  Brief mentions of renewal reward processes can also be found in Taylor and Karlin \cite{taylor-ISM} as well as in Heyman and Sobel \cite{heyman-SM} in a contribution by Serfozo. In Moder and Elmaghraby \cite{moder-HOR} one may find a short treatment of cumulative processes that also includes asymptotic results not only for the expectation of $R(t)$ as given by \eqref{eq:exp}, but also for its variance. As noted in the text there, there are multivariate versions of these results as well. The following books offer a thorough treatment of this topic and multitude of examples; moreover, they are approachable to a wide audience with only undergraduate knowledge of probability theory \cite{kulkarni-MASS,medhi-SP,resnick-ASP,ross-SP,tijms-SMA,wolff-SMTQ}.

\phantomsection
\addcontentsline{toc}{section}{References}

\end{document}